\documentclass[graybox]{svmult}

\usepackage{type1cm}        
\usepackage{makeidx}         
\usepackage{graphicx}        
\usepackage{multicol}        
\usepackage[bottom]{footmisc}
\usepackage{amsmath}

\usepackage{newtxtext}       %
\usepackage{newtxmath}       

\usepackage{color}

\usepackage{tikz}

\newcommand{\R}{\mathbb{R}}
\newenvironment{system}%
{\left\lbrace\begin{array}{r@{\hspace{1mm}}ll}}%
{\end{array}\right.}

\definecolor{Darkgblue}{rgb}{0.3,0.3,0.5}

\makeindex             

\begin{document}

\title*{On the unique continuation property of solutions to the two-dimensional Zakharov-Kuznetsov equation}
\titlerunning{Unique continuation property of ZK equation}
\author{Lucrezia Cossetti}
\institute{Lucrezia Cossetti \at Fakult\"at f\"ur Mathematik, Institut f\"ur Analysis, Karlsruher Institut f\"ur Technologie, Englerstra{\ss}e 2,76131 Karlsruhe, Germany, \email{lucrezia.cossetti@kit.edu}}
%
%
\maketitle


\abstract{The purpose of the current paper is twofold: to some extent it is intended as a \emph{review} of the recent optimal result in~\cite{C_F_L} concerning the unique continuation property of solutions to the two-dimensional Zakharov-Kuznetsov equation. On the other hand, the main core of the work is devoted to providing an \emph{alternative proof} of the aforementioned result. The importance of this original contribution relies on the fact that, unlike the approach used in~\cite{C_F_L}, the strategy adopted here is not sensitive of the two dimensional setting of the problem and therefore could be adapted to higher dimensional Zakharov-Kuznetsov equations for which, as far as we know, a proof of an analogous optimal unique continuation principle is still missing. For sake of clearness we focus here on the 2D case only, the higher dimensional analysis will be discussed somewhere else.}

\section{Introduction}

This paper is concerned with the study of unique continuation properties of solutions to the Zakharov-Kuznetsov (ZK) equation
\begin{equation}
\label{ZK}
	\partial_t u + \partial_x^3 u + \partial_x \partial_y^2 u + u \partial_x u=0, 
	\qquad (x,y)\in \R^2,
	t\in [0,1]. 
\end{equation} 
The equation was introduced in the context of plasma physic by Zakharov and Kuznetsov in~\cite{Z_K}, where they formally deduced that the propagation of nonlinear ion-acoustic waves in magnetized plasma is governed by this mathematical model. A rigorous derivation of equation~\eqref{ZK} was given by Lannes, Linares and Saut in~\cite{L_L_S}.

\medskip
The problem of local and global well-posedness for the Cauchy problem associated to~\eqref{ZK} has extensively been studied. Up to date the best local well-posedness result available in the literature was obtained independently by Molinet and Pilod~\cite{M_P} and Gr\"unrock and Herr~\cite{G_H} for initial data in $H^s(\R^2),$ $s>\frac{1}{2}.$
Then the global theory follows by standard arguments based on $L^2$ and $H^1$ conservation laws. We refer to~\cite{F, L_P, L_P_S, L_PII} and references therein for other results of this type and several additional remarks concerning with properties of this equation. 

\medskip
In this work we are mainly interested in the unique continuation setting of problems linked to the equation. At this regard we should mention the recent work~\cite{C_F_L}. The main achievement in~\cite{C_F_L} has been to obtain \emph{sharp} sufficient conditions on the (spatial) decay of the difference of two solutions of~\eqref{ZK} at two different times which guarantee that both solutions coincide. More precisely in~\cite{C_F_L} the following theorem is proved.
\begin{theorem}
\label{main}
	Suppose that for some small $\varepsilon>0$
	\begin{equation*}
		u_1, u_2 \in C([0,1]; H^4(\R^2)\cap L^2((1+ \left|x\right|)^{2(4/3 + \varepsilon)}dx dy)) \cap C^1([0,1]; L^2(\R^2)),
	\end{equation*}
	are solutions of~\emph{\eqref{ZK}}. Then there exists a universal constant $a_0>0$ such that if for some $a>a_0$
	\begin{equation}
	\label{decay}
		u_1(0)- u_2(0), u_1(1)-u_2(1)\in L^2(e^{a \left|x\right|^{3/2}} dx dy),
	\end{equation}
	then $u_1=u_2.$
\end{theorem} 
The exponential decay in condition~\eqref{decay} can be seen as a reflection of the asymptotic behavior of the fundamental solution to the associated linear operator $\partial_t + \partial_x^3 + \partial_x \partial_y^2,$ that is evidently given by the formula
\begin{equation*}
	G(t,x,y)= \theta(t) \mathcal{F}^{-1} [e^{i t(\xi^3 + \xi \eta^2)}]= \frac{\theta(t)}{t^{2/3}} S\left(\frac{x}{t^{1/3}}, \frac{y}{t^{1/3}}\right),
\end{equation*}
where 
\begin{equation*}
	S(x,y)= \frac{1}{2\pi} \mathcal{F}^{-1} [e^{i(\xi^3 + \xi \eta^2)}],
\end{equation*}
$\mathcal{F}^{-1}$ denotes the inverse Fourier transform in $\R^2,$ 
$\xi$ and $\eta$ are the variables in the frequency space corresponding to the space variables $x$ and $y$ respectively and $\theta$ is the Heaviside function. Indeed in~\cite{F_A} it has been proved  that $S$ displays an exponential decay just in the $x$ variable and more precisely it satisfies the estimate
\begin{equation*}
	\left|D^\nu S(x,y)\right|\leq c(1+ \left|y\right|)^{-m} e^{-c x_+^{3/2}}
	\qquad x_+:=\max\{x;0\}, y\in \R,
\end{equation*} 
for some constant $c>0$ and for any integer $m\geq 0$ and multi-index $\nu.$

Encouraged by this fact, in~\cite{C_F_L} the authors also show that the decay rate in Theorem~\ref{main} is optimal.

\medskip
The proof of Theorem~\ref{main} follows the scheme developed in~\cite{E_K_P_V} for the KdV equation and is based upon the comparison of two types of estimates, a lower bound, which comes as a consequence of a Carleman estimate for suitable compactly supported functions and that represents the most involved part of the proof, and an upper bound for the $L^2$-norm of the solution and its spatial derivatives up to order two, which exploits the exponential decay assumed in~\eqref{decay} for the data. 

The main obstacles that arose in trying to adapt the approach used for this one dimensional model, that is the KdV equation, to its higher dimensional generalization represented by the ZK equation, were basically due to the non symmetric form of this model in the two spatial variables $x, y.$ 

In order to overcome this difficulty the authors took advantage of the remarkable observation in~\cite{G_H} in which it was pointed out that after a suitable rotation of the set of coordinates, namely introducing the following change of variables
\begin{equation}
\label{change_variable}
	\begin{system}
		x'= \mu x + \lambda y\\
		y'= \mu x - \lambda y
	\end{system}
\end{equation} 
with $\mu = 4^{-1/3}$ and $\lambda= \sqrt{3}\mu,$ one is given with a \emph{symmetrized} version of the ZK equation~\eqref{ZK}. More precisely, defining $\widetilde{u}(x',y')=u(x,y),$ it is a simple computation to show that   
\begin{equation*}
	(\partial_x^3 + \partial_x \partial_y^2)u(x,y)= (\partial_{x'}^3 + \partial_{y'}^3) \widetilde{u}(x',y').
\end{equation*}
Thus one is led to consider the following equation
\begin{equation}
\label{ZKsymmetric}
	\partial_t u + (\partial_x^3 + \partial_y^3) u + 4^{-1/3} u (\partial_x + \partial_y)u=0, \qquad (x,y)\in \R^2, t\in [0,1].
\end{equation}
instead of~\eqref{ZK}. The main advantages of equation~\eqref{ZKsymmetric} over the classical form of~\eqref{ZK} are the aforementioned symmetric dependance upon $x$ and $y$ and the fact that the new equation comes with a dispersive part of KdV-type for both the spatial variables. These facts suggest that the expected optimal exponential decay rate to get uniqueness for equation~\eqref{ZKsymmetric} should be symmetric in $x$ and $y$ and with as exponent the one inherited from the asymptotic decay of the Airy function, i.e. $3/2.$ This claim found support and confirmation in~\cite{C_F_L} where the following unique continuation principle was obtained for the symmetric ZK equation~\eqref{ZKsymmetric}. 
\begin{theorem}
\label{thm:intermediate}
	Suppose that for some small $\varepsilon>0$
	\begin{equation*}
		u_1, u_2 \in C([0,1]; H^4(\R^2)\cap L^2((1+ \left|x + y\right|)^{2(4/3 + \varepsilon)}dx dy)) \cap C^1([0,1]; L^2(\R^2)),
	\end{equation*}
	are solutions of~\emph{\eqref{ZKsymmetric}}. Then there exists a universal constant $a_0>0$ such that if for some $a>a_0$
	\begin{equation}
	\label{decay_symmetric}
		u_1(0)- u_2(0), u_1(1)-u_2(1)\in L^2(e^{a \left|x+y\right|^{3/2}} dx dy),
	\end{equation}
	then $u_1=u_2.$
\end{theorem} 

Observe that an alternative to the decay rate assumed in~\eqref{decay_symmetric}, that is still symmetric in $x$ and $y$ and that enjoys the same $3/2$-homogeneity, would have been requiring 
\begin{equation}
\label{alternative}
	u_1(0)- u_2(0), u_1(1)-u_2(1)\in L^2(e^{a (x^2+y^2)^{3/4}} dx dy).
\end{equation} 
Nevertheless one realizes that this choice does not work as well as~\eqref{decay_symmetric}. In order to see that, first we stress that, even if Theorem~\ref{thm:intermediate} is of its independent interest, the original purpose in~\cite{C_F_L} was to pass through the proof of Theorem~\ref{thm:intermediate} only as an intermediate step to get Theorem~\ref{main}. In other words, the final goal in~\cite{C_F_L} was to deduce Theorem~\ref{main} from Theorem~\ref{thm:intermediate} once one gets back to the original variables. In particular this requires that hypothesis~\eqref{decay_symmetric} (or alternatively~\eqref{alternative} ``translates'' well to~\eqref{decay} after the inverse change of variables. It is precisely at this stage that the alternative assumption~\eqref{alternative} fails to work. In order to see that we define
\begin{equation*}
	\widetilde{\varphi}(x',y'):=|x'+ y'|.
\end{equation*}
Coming back to the starting set of coordinates, namely using~\eqref{change_variable}, we have
\begin{equation*}
	\widetilde{\varphi}(x',y'):=|x'+ y'|=2\mu |x|=2\mu \varphi(x),
\end{equation*}
with $\varphi(x):=|x|.$ Clearly one observes that the form of $\varphi$ is consistent with the exponential decay (\emph{only} in $x$) contained in assumption~\eqref{decay} of Theorem~\ref{main}. 
On the contrary, defining
\begin{equation*}
	\psi(x',y'):=x'^2+ y'^2,
\end{equation*}
it is easy to see that $\psi$ is almost invariant under the change of variable~\eqref{change_variable}, more precisely
\begin{equation*}
	\psi(x', y'):=x'^2 + y'^2
	=2\mu^2 x^2 +2\lambda^2 y^2=\psi(\sqrt{2}\mu x, \sqrt{2}\lambda y),
\end{equation*}
from this one realizes that there is no chance to get from the alternative choice~\eqref{alternative} the corresponding desired and optimal decay in~\eqref{decay}. 

In this respect we should mention that Theorem~\ref{main} improved the following previous result contained in~\cite{B_I_M}.
\begin{theorem}
\label{BIM}
	Suppose that for some small $\varepsilon>0$
	\begin{equation*}
		u_1,u_2 \in C([0,1]; H^4(\R^2)\cap L^2((1+ x^2 + y^2)^{4/3 + \varepsilon}dx dy)) \cap C^1([0,1]; L^2(\R^2)),
	\end{equation*}
	are solutions of~\emph{\eqref{ZK}}. Then there exists a universal constant $a_0>0$ such that if for some $a>a_0$
	\begin{equation}
	\label{decayBIM}
		u_1(0)- u_2(0), u_1(1)-u_2(1)\in L^2(e^{a (x^2+y^2)^{3/4}} dx dy),
	\end{equation}
	then $u_1=u_2.$
\end{theorem} 
From assumption~\eqref{decayBIM} it is clear that it is precisely this result the one that we would have obtained proving Theorem~\ref{thm:intermediate} for the symmetric ZK equation~\eqref{ZKsymmetric} with the alternative hypothesis~\eqref{alternative} and then coming back to the original equation~\eqref{ZK} with the inverse change of variables. In other words, the authors in~\cite{B_I_M} guessed the right $3/2$-homogeneity of the optimal decay rate, but they missed the correct form of the exponential weight.
    
\medskip
In passing we stress that avoiding to assume exponential decay of the data in the $y$ variable (see condition~\eqref{decay} produced new highly non trivial technical difficulties in the proof of Theorem~\ref{main} compared to Theorem~\ref{BIM}. One of the main difficulty arises in the use of the general Carleman estimate to obtaining the mentioned lower bound for our solution of~\eqref{ZK}. As expected if one wants to apply such a lower bound to prove a unique continuation property, this requires the introduction of suitable cut-off functions whose supports have to be ``related'' to the sublevel sets of the exponent-functions in the decay assumptions for the data (respectively $\varphi(x)=|x|$ in Theorem~\ref{main} and $\psi(x,y)=x^2+y^2$ in Theorem~\ref{BIM}). Now, it should appear clear that due to the \emph{unboundedness} in $\R^2$ of the sublevel set associated to $\varphi,$ in contrast with the ball-type set associated to $\psi,$ a more careful analysis is needed in this case in order to meaningfully apply the Carleman estimate valid for compactly supported functions only.
At this proposal, for sake of clearness we will take advantage of the Appendix in this manuscript to fix a small glitch we found in~\cite{C_F_L} although this is not affecting sensitively the proof of the main result there.
   
\bigskip
Recently there has been a growing interest in studying higher dimensional Zakharov-Kuznetsov equations. In particular a strong impetus has been given to the investigation from diverse perspectives of the $(1+ 3)$-dimensional ZK which has the following form
\begin{equation}
\label{3dZK}
	\partial_t u + \partial_x \Delta u + u \partial_x u=0, 
	\qquad (x,y,z)\in\R^3, t\in [0,1], 
\end{equation}
where $\Delta$ is the three dimensional spacial Laplace operator. 
Among the different aspects in the study, the proof of a unique continuation principles for such an evolution equation has consistently attracted the community. As far as we know the most up-to-date result in this direction is contained in~\cite{B_J_M} and is stated as follows:
\begin{theorem}
\label{3d}
	Suppose for some small $\varepsilon>0$
	\begin{equation*}
		u_1,u_2\in C([0,1]; H^4(\R^3)\cap L^2((1+ x^2 + y^2 + z^2)^{8/5 + \varepsilon} dx dy dz)) \cap C^1([0,1]; L^2(\R^3)),
	\end{equation*}
	are solutions to~\emph{\eqref{3dZK}}. Then there exists a universal constant $a_0>0,$ such that if for some $a>a_0$ 
	\begin{equation}
	\label{decay3d}
		u_1(0)-u_2(0), u_1(1)-u_2(1)\in L^2(e^{a(x^2+y^2+z^2)^{3/4}} dx dy dz),
	\end{equation}
	then $u_1=u_2.$
\end{theorem}
As one can immediately notice, Theorem~\ref{3d} is the natural generalization of Theorem~\ref{BIM} (compare~\eqref{decay3d} with~\eqref{decayBIM}) and, as Theorem~\ref{BIM} itself, cannot be optimal due to the symmetric assumption~\eqref{decay3d} that does not reflect the non symmetric form of equation~\eqref{3dZK}. This means that for this three dimensional model a proof of a unique continuation principle with \emph{sharp} decay assumptions on the data is still missing. We also stress that the change of coordinates~\eqref{change_variable} introduced in~\cite{C_F_L} which led to the optimal result in the two dimensional setting does not generalize to the 3D case. This motivates us in trying to re-prove Theorem~\ref{main} avoiding the change of variables which deviated the study to the symmetric ZK equation~\eqref{ZKsymmetric} and working, instead, on the 
original ZK equation~\eqref{ZK} itself. This alternative proof represents the main part of these notes and is contained in the next section. For sake of clearness we will focus here on the 2D case only, the higher dimensional generalization of the proof provided here will appear somewhere else.

\section{Proof of Theorem~\ref{main}}
This section is concerned with the proof of Theorem~\ref{main}. As customarily in proving unique continuation results, this will follow as a consequence of a comparison between suitable lower and upper bounds for the difference $v:=u_1-u_2$ of two solutions $u_1, u_2$ of~\eqref{ZK}. Working with the difference $v$ addresses the problem of proving uniqueness for the equation~\eqref{ZK} into the problem of proving the triviality of the solution to the corresponding equation associated to $v,$ namely
\begin{equation}
\label{ZKv}
\partial_t v + \partial_x^3 v + \partial_x \partial_y^2 v + u_1 \partial_x v + \partial_x u_2 v=0.  
\end{equation}
Actually for our analysis we consider a more general problem than~\eqref{ZKv}, that is
\begin{equation}
\label{ZKv_general}
	\partial_t v + \partial_x^3 v + \partial_x \partial_y^2 v + a_1(x,y,t) \partial_x v + a_0(x,y,t) v=0,
\end{equation}
for suitable $a_0, a_1.$ Clearly~\eqref{ZKv} is a particular case of~\eqref{ZKv_general} setting $a_0=\partial_x u_2$ and $a_1=u_1.$
 
We start introducing explicitly the quantity that we aim at estimating from below and from above, namely 
\begin{equation}
\label{ARv}
	A_R(v):=\left(\int_0^1 \int_{Q_R} (|v|^2 + |\partial_x v|^2 + |\partial_y v|^2 + |\partial_x \partial_y v|^2 + |\Delta v|^2) dx dy dt \right)^{1/2},
\end{equation}
where $Q_R:= \{(x, y)\colon R-1\leq |x|\leq R \wedge |y|\leq2\}.$

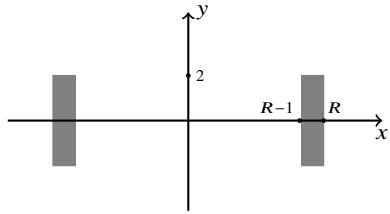
\begin{figure}[ht]
\centering
\begin{tikzpicture}[scale=.3]
\filldraw[gray] (5,2)--(6,2)--(6,-2)--(5,-2);
\filldraw[gray] (-5,2)--(-6,2)--(-6,-2)--(-5,-2);
\draw[domain=5:6, color=gray, samples=2] plot(\x, {2});
\coordinate [label=above left: \textcolor{black}{$\scriptstyle R-1$}] (a) at (5,0); 
\fill[black] (4.93,0) circle (3pt);
\coordinate [label=above right: \textcolor{black}{$\scriptstyle R$}] (b) at (5.8,0); 
\fill[black] (6,0) circle (3pt);
\coordinate [label=right: \textcolor{black}{$\scriptstyle 2$}] (c) at (0,2); 
\fill[black] (0,2) circle (3pt);
\draw[thick,->](-8,0)--(8.6,0) node[below]{$x$};
\draw[thick,->](0,-4)--(0,4.8) node[right]{$y$};
\end{tikzpicture}
\caption{The region $Q_R$} 
\label{fig:an_dom}
\end{figure}

The statements of the lower and upper bound are contained in the following two theorems.
\begin{theorem}[Lower bound]
\label{lower_bound}
	Let $v\in C([0,1]; H^3(\R^2))\cap C^1([0,1];L^2(\R^2))$ be a solution to~\eqref{ZKv_general} with $a_0, a_1\in L^\infty(\R^3).$ Assume that there exists a positive constant $A$ such that 
\begin{equation*}	
	\int_0^1 \int_{\R^2} (|v|^2 + |\partial_x v|^2 + |\partial_y v|^2 + |\Delta v|^2 )dxdydt \leq A^2.
\end{equation*}
	
	Let $\delta>0,$ $r\in (0,1/2)$ and $Q:=\{(x,y,t)\colon \sqrt{x^2+y^2}\leq 1, t \in [r,1-r]\}$ and suppose that $||v||_{L^2(Q)}>\delta.$ 
	
	Then there exist positive constants $\widetilde{R}_0, c_0, c_1$ depending on $A, ||a_0||_{L^\infty(\R^3)}$ and $||a_1||_{L^\infty(\R^3)}$ such that for $R\geq \widetilde{R}_0$
	\begin{equation}
	\label{eq:lower_bound}
		A_{R}(v)\geq c_0e^{-c_1 R^{3/2}}.
	\end{equation}
\end{theorem}
\begin{theorem}[Upper bound]
\label{upper_bound}
Let $v\in C([0,1]; H^4(\R^2))$ be a solution of~\eqref{ZKv_general} whose coefficients $a_0$ and $a_1$ satisfy $a_0\in L^\infty \cap L^2_x L^\infty_{yt}$ and $a_1\in L^2_x L^\infty_{yt} \cap L^1_x L^\infty_{yt}$ respectively. Then there exist positive constants $c$ and $R_0,$ such that if for some $a>0$
\begin{equation*}
	v(0), v(1)\in L^2(e^{a|x|^{3/2}}dx dy),
\end{equation*} 
then 
\begin{equation}
\label{eq:upper_bound}
	A_R(v)\leq c e^{-a \big(\frac{R}{18}\big)^{3/2}}
\end{equation}
holds true for $R\geq R_0.$
\end{theorem}

Before proving Theorem~\ref{lower_bound} and Theorem~\ref{upper_bound} we show how Theorem~\ref{main} follows easily as a consequence of these two deep results.

\begin{proof}(of Theorem~\ref{main})
	We argue by contradiction. Suppose that $v:=u_1-u_2\neq 0.$ Then we can assume, after a possible translation, dilation and multiplication by a constant, that $v$ satisfies the hypotheses of Theorem~\ref{lower_bound}. Moreover, it can be easily seen (see~\cite{C_F_L} for details) that fixing $a_0=\partial_x u_2$ and $a_1=u_1,$ then $a_0\in L^\infty \cap L^2_x L^\infty_{yt}$ and $a_1\in L^2_x L^\infty_{yt} \cap L^1_x L^\infty_{yt}$ and so the hypotheses of Theorem~\ref{upper_bound} are satisfied too. Thus, combining~\eqref{eq:lower_bound} and~\eqref{eq:upper_bound}, one has for sufficiently large $R$
	\begin{equation*}
		c_0e^{-c_1 R^{3/2}}\leq A_R(v)\leq c e^{-\frac{a}{18^{3/2}}R^{3/2}}.
	\end{equation*}
			Finally, assuming $a>a_0:=18^{3/2} c_1$ and taking the limit as $R$ tends to infinity, we get a contradiction. Therefore $v=0$ and Theorem~\ref{main} is proved. 
\end{proof}

\subsection{Lower bound: Proof of Theorem~\ref{lower_bound}}
The starting point in the proof of the lower bound in Theorem~\ref{lower_bound} is a Carleman estimate for the operator
\begin{equation*}
	P=\partial_t + \partial_x^3 + \partial_x\partial_y^2 + a_1(x,y,t)\partial_x + a_0(x,y,t),
\end{equation*} 
with $a_0, a_1\in L^\infty(\R^3).$ At this proposal we take advantage of the corresponding Carleman estimate for the leading part of $P,$ namely $\partial_t + \partial_x^3 + \partial_x \partial_y^2,$ already proved in~\cite{B_I_M}. We restate here the result for the reader convenience. 
\begin{lemma}
	Let $\varphi\colon [0,1]\to \R$ be a smooth function and let $D:=\R^2\times [0,1].$ Let us assume that $R>1$ and define
	\begin{equation*}
		\phi(x,y,t):= \left(\frac{x}{R} + \varphi(t) \right)^2 + \frac{y^2}{R^2}.
	\end{equation*}
	Then there exists $\overline{C}= \max\{|| \varphi'||_{L^\infty}, || \varphi''||_{L^\infty}, 1\}>0$ such that the inequality
	\begin{equation}
	\label{Carleman_leading}
		\frac{\alpha^{5/2}}{R^3}||e^{\alpha \phi} g||_{L^2(D)} + \frac{\alpha^{3/2}}{R^2}||e^{\alpha \phi} \partial_x g||_{L^2(D)}
		\leq \sqrt{2} ||e^{\alpha \phi} (\partial_t + \partial_x^3 + \partial_x \partial_y^2) g||_{L^2(D)}
	\end{equation}
	holds if $\alpha\geq \overline{C}R^{3/2}$ and $g\in C([0,1]; H^3(\R^2)) \cap C^1([0,1];L^2(\R^2))$ is compactly supported in $\{(x,y,t)\colon | \frac{x}{R} + \varphi(t) |\geq 1\}.$ 
\end{lemma}

From estimate~\eqref{Carleman_leading} the corresponding inequality for the full operator $P$ is readily given. Indeed, adding and subtracting the lower order terms and using the H\"older inequality one has
\begin{multline}
		\frac{\alpha^{5/2}}{R^3}||e^{\alpha \phi} g||_{L^2(D)} + \frac{\alpha^{3/2}}{R^2}||e^{\alpha \phi} \partial_x g||_{L^2(D)}\\
		\leq \sqrt{2} ||e^{\alpha \phi} (\partial_t + \partial_x^3 + \partial_x\partial_y^2 + a_1\partial_x + a_0)g||_{L^2(D)}
		\\ + \sqrt{2} ||a_1||_{L^\infty(D)} ||e^{\alpha \phi} \partial_x g||_{L^2(D)}
		\\ + \sqrt{2} ||a_0||_{L^\infty(D)} ||e^{\alpha \phi} g||_{L^2(D)}.
\end{multline} 
Observe that under the hypothesis $\alpha \geq \overline{C} R^{3/2},$ the ratios $\alpha^{3/2}/R^2$ and $\alpha^{5/2}/R^3$ grow as a positive fractional power of $R,$ as a consequence,  as soon as $R$ is taken sufficiently large, 
the last two terms on the right-hand side can be absorbed into the corresponding terms in the left-hand side giving the desired estimate. 

Summing up, we have proved the following result:
\begin{lemma}
\label{lemma:Carleman}
	Let $\varphi\colon [0,1]\to \R$ be a smooth function and let $D:=\R^2\times [0,1].$ Let us define
	\begin{equation*}
		\phi(x,y,t):= \left(\frac{x}{R} + \varphi(t) \right)^2 + \frac{y^2}{R^2}.
	\end{equation*}
	Then there exist constants $c>0,$ $R_0=R_0(||a_0||_{L^\infty}, ||a_1||_{L^\infty})>1$ and  $\overline{C}= \max\{|| \varphi'||_{L^\infty}, || \varphi''||_{L^\infty}, 1\}>0$ such that the inequality
	\begin{multline}
	\label{Carleman}
		\frac{\alpha^{5/2}}{R^3}||e^{\alpha \phi} g||_{L^2(D)} + \frac{\alpha^{3/2}}{R^2}||e^{\alpha \phi} \partial_x g||_{L^2(D)}
		\\\leq c||e^{\alpha \phi} (\partial_t + \partial_x^3 + \partial_x \partial_y^2 + a_1 \partial_x + a_0) g||_{L^2(D)}
	\end{multline}
	holds if $\alpha\geq \overline{C}R^{3/2},$ $R\geq R_0$ and $g\in C([0,1]; H^3(\R^2)) \cap C^1([0,1];L^2(\R^2))$ is compactly supported in $\{(x,y,t)\colon | \frac{x}{R} + \varphi(t) |\geq 1\}.$ 
\end{lemma}

Now we are in position to prove the lower bound in Theorem~\ref{lower_bound}.
\begin{proof}[of Theorem~\ref{lower_bound}]
Let us introduce the function $\theta_{R,\varepsilon}\in C^\infty_0(\R^2)$ such that $\theta_{R,\varepsilon}(x,y)=1$ on $\{(x,y)\colon |x|< R-1 \wedge |y|< 2\}$ and $\theta_{R,\varepsilon}(x,y)=0$ on $\R^2\setminus \{(x,y)\colon |x|\leq R \wedge |y|\leq 2+ \varepsilon\}.$ Let $\mu \in C^\infty(\R)$ be such that $\mu(x)=0$ if $|x|< 1$ and $\mu(x)=1$ if $|x|> 2,$ and $\varphi\colon \R \to [0,2\sqrt{2}]$ such that $\varphi(t)=0$ on $[0, r/2] \cup [1-r/2,1],$ $\varphi(t)=2\sqrt{2}$ on $[r,1-r],$ increasing in $[r/2, r]$ and decreasing in $[1-r,1-r/2].$

We define 
\begin{equation}
\label{g}
	g(x,y,t)=\theta_{R, \varepsilon}(x, y)\mu \left(\frac{x}{R} + \varphi(t) \right) v(x,y,t), \qquad (x,y)\in \R^2, t\in [0,1],
\end{equation}
where $v$ satisfies~\eqref{ZKv_general}.

It is easy to see that $g$ satisfies
\begin{equation}
\label{eq_g}
	(\partial_t + \partial_x^3 + \partial_x \partial_y^2 + a_1\partial_x + a_0)g=F_1 + F_2,
\end{equation}
where
\begin{equation*}
	\begin{split}
F_1= \mu\left(\frac{x}{R} + \varphi(t) \right) 
			\Big[ 
			&3\partial_x \theta_{R, \varepsilon} \partial_x^2 v 
			+2\partial_y \theta_{R, \varepsilon} \partial_x\partial_y v
			+\partial_x \theta_{R, \varepsilon} \partial_y^2 v\\
			&+3\partial_x^2 \theta_{R, \varepsilon} \partial_x v
			+\partial_y^2 \theta_{R, \varepsilon} \partial_x v
			+2\partial_x\partial_y \theta_{R, \varepsilon} \partial_y v\\
			&+\partial_x^3 \theta_{R, \varepsilon} v
			+\partial_x \partial_y^2 \theta_{R, \varepsilon} v
			+a_1 \partial_x \theta_{R, \varepsilon} v
			\Big]
	\end{split}
\end{equation*}
and 
\begin{equation*}
	\begin{split}
	F_2= &3 R^{-1}\theta_{R, \varepsilon} \partial_x \mu \partial_x^2 v
	+ R^{-1}\theta_{R, \varepsilon} \partial_x \mu \partial_y^2 v\\
	&+ R^{-1} (3 R^{-1}\theta_{R, \varepsilon} \partial_x^2 \mu + 6 \partial_x \theta_{R, \varepsilon} \partial_x \mu)\partial_x v
	+ 2 R^{-1}\partial_y \theta_{R, \varepsilon} \partial_x \mu \partial_y v\\
	&+ \Big[ 
	\theta_{R, \varepsilon} \partial_x \mu (\varphi' + a_1 R^{-1})
	+ R^{-1} (3 \partial_x^2 \theta_{R, \varepsilon} +  \partial_y^2 \theta_{R, \varepsilon}) \partial_x \mu\\
	&\phantom{+\Big[}+ 3 R^{-2} \partial_x \theta_{R, \varepsilon} \partial_x^2 \mu
	+ \theta_{R, \varepsilon} \partial_x^3 \mu 
	\Big]v.
	\end{split}
\end{equation*}
First we show that $g$ as defined in~\eqref{g} satisfies the hypotheses of Lemma~\ref{lemma:Carleman}.
\begin{itemize}
	\item if $(x, y)$ belongs to $\R^2\setminus\{(x,y)\colon |x|\leq R \wedge |y|\leq 2+ \varepsilon\}$ we are outside the support of $\theta_{R, \varepsilon}$ and then $g=0.$
	\item if $(x, y)$ belongs to $\{(x,y)\colon |x|< R \wedge |y|< 2+ \varepsilon\}$ and $t\in [0,r/2]\cup[1-r/2,1]$ then $g=0.$ Indeed being $|x|< R$ and since $\varphi(t)=0$ if $t\in [0,r/2]\cup[1-r/2,1],$ then $|\frac{x}{R} + \varphi(t)|<1,$ therefore we are out of the support of $\mu(\frac{x}{R} + \varphi(t))$ and so $g=0.$
\end{itemize}
This guarantees that $g$ is compactly supported. Now we observe that $g$ is supported in $\{(x,y,t)\colon |\frac{x}{R}+ \varphi(t)|\geq 1\}.$ Indeed if $|\frac{x}{R}+ \varphi(t)|< 1$ then $\mu(\frac{x}{R} + \varphi(t))=0$ and so $g=0.$

Applying Carleman estimate~\eqref{Carleman} of Lemma~\ref{lemma:Carleman} to $g$ and using~\eqref{eq_g} gives
\begin{equation}
\label{intermediate}
	\begin{split}
	c \frac{\alpha^{5/2}}{R^3} ||e^{\alpha \phi} g||_{L^2(D)}
	&\leq ||e^{\alpha \phi}(\partial_t + \partial_x^3 + \partial_x \partial_y^2 + a_1 \partial_x + a_0) g||_{L^2(D)}\\
	&\leq ||e^{\alpha \phi} F_1||_{L^2(D)} + ||e^{\alpha \phi} F_2||_{L^2(D)}.
	\end{split}
\end{equation} 
Now we comment further on the terms involving $F_1$ and $F_2.$ As regards with $F_1,$ one observes that all the terms in $F_1$ contain derivatives of $\theta_{R, \varepsilon},$ from the definition of $\theta_{R, \varepsilon}$ this implies that $F_1$ is supported in the set $Q_{R, \varepsilon}\times [0,1],$ where
$Q_{R, \varepsilon}$ is a ``anulus''-type domain defined as follows:
\begin{equation*}
	Q_{R,\varepsilon}
	:=\{(x,y)\colon |x|\leq R \wedge |y|\leq 2+ \varepsilon\} 
	\setminus 
	\{(x,y)\colon |x|< R-1 \wedge |y|<2	\}.
\end{equation*} 
It is not difficult to show that in $Q_{R, \varepsilon}\times [0,1]$ one has $\phi(x,y,t)\leq 20.$ In order to see that observe that if $(x,y)\in Q_{R, \varepsilon}$ then $\sqrt{x^2+y^2}\leq \sqrt{2}R.$ As a consequence, using the trivial inequality $2ab\leq a^2 + b^2,$ one has
\begin{equation*}
	\phi(x,y,t):= \left(\frac{x}{R} + \varphi(t) \right)^2 + \frac{y^2}{R^2}\leq 2\left(\frac{x^2}{R^2} + \frac{y^2}{R^2} \right) + 2\varphi(t)^2\leq 20.
\end{equation*}
With respect to $F_2,$ since the derivatives of $\mu$ appear, $F_2$ is supported in $\{(x,y,t)\colon 1\leq|\frac{x}{R} + \varphi(t)|\leq 2, t\in [0,1]\}.$ Moreover, since we are still in the support of $\theta_{R, \varepsilon},$ in particular $\sqrt{x^2 + y^2} \leq \sqrt{2} R.$ From this, it follows that $\phi(x,y,t)\leq 6,$ indeed
\begin{equation*}
	\phi(x,y,t):=\left(\frac{x}{R} + \varphi(t) \right)^2 + \frac{y^2}{R^2} \leq 4+ 2= 6.
\end{equation*}
Using the last remarks, from inequality~\eqref{intermediate} one gets  
\begin{equation}
\label{first}
c \frac{\alpha^{5/2}}{R^3} ||e^{\alpha \phi} g||_{L^2(D)}
\leq c_1 e^{20 \alpha} A_{R, \varepsilon}(v) 
+ c_2 e^{6\alpha} A,
\end{equation}
where we have defined
\begin{equation*}
	A_{R, \varepsilon}(v)= \left( \int_0^1 \int_{Q_{R, \varepsilon}} (|v|^2 + |\partial_x v|^2 + |\partial_y v|^2 + |\partial_x \partial_y v|^2 + |\Delta v|^2) dx dy dt  \right)^{1/2}.	
\end{equation*}
Observe that in $Q:=\{(x,y,t)\colon \sqrt{x^2+y^2}\leq 1, t\in [r,1-r]\}$ we have $g=v.$ Indeed if $\sqrt{x^2+y^2}\leq 1$ and $R$ is sufficiently large, then $\theta_{R, \varepsilon}=1.$ Moreover, in $Q$ it also holds true that $\mu(\frac{x}{R} + \varphi(t))=1.$ Indeed, using that $\phi(t)=2\sqrt{2}$ in $[r, 1-r]$ and the trivial inequality $a+b\leq \sqrt{2}\sqrt{a^2 + b^2}$ with $a,b\geq 0,$ if $R$ is sufficiently large, say $R\geq 4,$ one has
\begin{equation}
\label{mu1}
	\left(\frac{x}{R} + \phi(t) \right)^2= \frac{x^2}{R^2} + 8 + \frac{4\sqrt{2}}{R}x\geq 8 -\frac{8}{R} \sqrt{x^2 + y^2}\geq 6> 4.
\end{equation}  
Then $|\frac{x}{R} + \varphi(t)|>2$ and so $\mu(\frac{x}{R} + \varphi(t))=1.$

Using that $g=v$ in $Q$ we obtain the following chain of inequalities:
\begin{equation}
\label{second}
	\begin{split}
	c \frac{\alpha^{5/2}}{R^3} ||e^{\alpha \phi} g||_{L^2(D)}
	&\geq c \frac{\alpha^{5/2}}{R^3} ||e^{\alpha \phi} g||_{L^2(Q)}
	=c \frac{\alpha^{5/2}}{R^3} ||e^{\alpha \phi} v||_{L^2(Q)}\\
	&\geq 
	c \frac{\alpha^{5/2}}{R^3} e^{6\alpha} ||v||_{L^2(Q)},
	\end{split}
\end{equation}
where in the last inequality we have used that in $Q$ one has $\phi\geq 6$ which follows reasoning as in~\eqref{mu1}. 

Using~\eqref{first} and~\eqref{second}, the assumption $||v||_{L^2(Q)}>\delta$ and dividing by $e^{6\alpha}$ one gets
\begin{equation*}
	c \frac{\alpha^{5/2}}{R^3} \delta\leq c_1 e^{14\alpha} A_{R, \varepsilon}(v) + c_2 A. 
\end{equation*} 
Taking $\alpha= \overline{C} R^{3/2}$ with $\overline{C}$ as in Lemma~\ref{lemma:Carleman} we obtain
\begin{equation*}
	c \overline{C}^{5/2} R^{3/4}\delta \leq c_1 e^{14 \overline{C}R^{3/2}} A_{R, \varepsilon}(v) + c_2 A. 
\end{equation*}
Now if we take $R$ large enough, the second term on the right-hand side can be absorbed into the left hand side. More precisely there exists $\widetilde{R}_0\geq R_0$ such that for $R\geq \widetilde{R}_0$ the following lower bound
\begin{equation}
\label{epsLB}
A_{R, \varepsilon}(v)\geq C e^{-14 \overline{C}R^{3/2}}
\end{equation}
holds true.
Now we write
\begin{equation*}
	Q_{R, \varepsilon}= Q_{R,\varepsilon}^{(1)} \cup Q_{R,\varepsilon}^{(2)},
\end{equation*}
where 
$
	Q_{R, \varepsilon}^{(1)}:= 
	\{(x, y)\colon |x|\leq R-1 \wedge 2\leq|y|\leq 2+ \varepsilon\}
$
	and
$
	Q_{R, \varepsilon}^{(2)}:= 
	\{(x, y)\colon R-1\leq|x|\leq R \wedge |y|\leq 2+ \varepsilon\}.
$

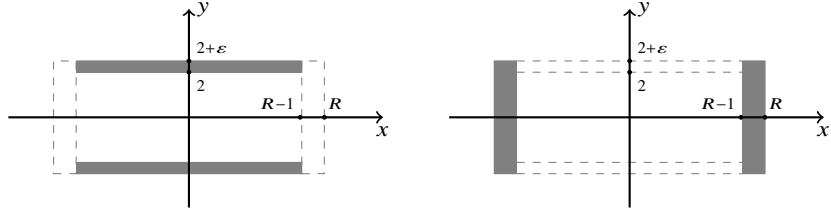
\begin{figure}[ht]
\centering
\begin{tikzpicture}[scale=.3]
\draw[fill=white, draw=gray, dashed] (5,-2.5) rectangle (6,2.5);
\draw[fill=white, draw=gray, dashed] (-6,-2.5) rectangle (-5,2.5);
\filldraw[gray] (-5,2.5)--(5,2.5)--(5,2)--(-5,2);
\filldraw[gray] (-5,-2.5)--(5,-2.5)--(5,-2)--(-5,-2);
\coordinate [label=above left: \textcolor{black}{$\scriptstyle R-1$}] (a) at (5,0); 
\fill[black] (4.93,0) circle (3pt);
\coordinate [label=above right: \textcolor{black}{$\scriptstyle R$}] (b) at (5.8,0); 
\fill[black] (6,0) circle (3pt);
\coordinate [label=below right: \textcolor{black}{$\scriptstyle 2$}] (c) at (0,2); 
\fill[black] (0,2) circle (3pt);
\coordinate [label=above right: \textcolor{black}{$\scriptstyle 2+ \varepsilon$}] (d) at (0,2.5); 
\fill[black] (0,2.5) circle (3pt);
\draw[thick,->](-8,0)--(8.6,0) node[below]{$x$};
\draw[thick,->](0,-4)--(0,4.8) node[right]{$y$};

\end{tikzpicture}
\qquad
\begin{tikzpicture}[scale=.3]
\filldraw[gray] (5,2.5)--(6,2.5)--(6,-2.5)--(5,-2.5);
\filldraw[gray] (-5,2.5)--(-6,2.5)--(-6,-2.5)--(-5,-2.5);
\draw[domain=-5:5, color=gray, dashed, samples=2] plot(\x, {2});
\draw[domain=-5:5, color=gray, dashed, samples=2] plot(\x, {2.5});
\draw[domain=-5:5, color=gray, dashed, samples=2] plot(\x, {-2});
\draw[domain=-5:5, color=gray, dashed, samples=2] plot(\x, {-2.5});
\coordinate [label=above left: \textcolor{black}{$\scriptstyle R-1$}] (a) at (5,0); 
\fill[black] (4.93,0) circle (3pt);
\coordinate [label=above right: \textcolor{black}{$\scriptstyle R$}] (b) at (5.8,0); 
\fill[black] (6,0) circle (3pt);
\coordinate [label=below right: \textcolor{black}{$\scriptstyle 2$}] (c) at (0,2); 
\fill[black] (0,2) circle (3pt);
\coordinate [label=above right: \textcolor{black}{$\scriptstyle 2+ \varepsilon$}] (d) at (0,2.5); 
\fill[black] (0,2.5) circle (3pt);
\draw[thick,->](-8,0)--(8.6,0) node[below]{$x$};
\draw[thick,->](0,-4)--(0,4.8) node[right]{$y$};

\end{tikzpicture}

\caption{The region $Q_{R,\varepsilon}^{(1)}$ (left) and  $Q_{R,\varepsilon}^{(2)}$ (right)} 
\end{figure}

It is easy to show that the measure of $Q_{R,\varepsilon}^{(1)}$ tends to zero as $\varepsilon$ goes to zero and that 
\begin{equation*}
	Q_{R, \varepsilon}^{(2)} \xrightarrow{\varepsilon \to 0} Q_{R},
\end{equation*}
where $Q_R:=\{(x, y)\colon R-1\leq|x|\leq R \wedge |y|\leq 2\}.$ 

It follows from the previous facts that, using the dominated convergence theorem, in the limit $\varepsilon \to 0$ one has  $A_{R, \varepsilon}(v) \to A_{R}(v),$ with $A_R(v)$ defined as in~\eqref{ARv}. 

Therefore, passing to the limit $\varepsilon \to 0$ in~\eqref{epsLB}, gives 
\begin{equation*}
A_R(v)\geq C e^{-14 \overline{C}R^{3/2}},
\end{equation*}
which is the desired lower bound.
\end{proof}

\subsection{Upper bound: Proof of Theorem~\ref{upper_bound}}
This section is concerned with the proof of the upper bound in Theorem~\ref{upper_bound}. The starting point is providing the following persistence-decay result, which, as we shall see, will come as consequence of boundedness properties for the inverse of the operator $P=\partial_t + \partial_x^3 + \partial_y\partial_x^2 + a_1\partial_x + a_0.$ 
\begin{lemma}
\label{lemma:persistence}
	Let $w\in C([0,1]; H^4(\R^2))\cap C^1([0,1];L^2(\R^2))$ such that for all $t\in [0,1],$ the support of $w(t)$ is contained in a compact subset $K$ of $\R^2$ and let $D:=\R^2\times [0,1].$
		
		Assume that $a_0\in L^\infty \cap L^2_xL^\infty_{y t}$ and $a_1\in L^2_x L^\infty_{yt}\cap L^1_xL^\infty_{yt},$ with small norms in these spaces.
		
		Then there exists $c>0$ independent of the set $K$ such that for $\lambda\geq 1$ the following estimate 
		\begin{multline}
		\label{eq:persistence}
			||e^{\lambda |x|}w||_{L^2(D)}
		+\sum_{0<k+l\leq 2}||e^{\lambda |x|}\partial_x^k \partial_y^l w||_{L^\infty_x L^2_{y t} (D)}\\
		\leq 
		c \lambda^2 \left(
		||J^3(e^{\lambda |x|} w(0)) ||_{L^2(\R^2)} 
		+ ||J^3(e^{\lambda |x|} w(1)) ||_{L^2(\R^2)}
		\right)\\
		+ c ||e^{\lambda |x|} (\partial_t + \partial_x^3 + \partial_y\partial_x^2 + a_1\partial_x + a_0) w ||_{L^1_t L^2_{xy}(D)\cap L^1_x L^2_{ty}(D)}
		\end{multline}
		holds true. 
		Here $J$ is defined in the Fourier space as $\mathcal{F}(Jg):= (1+ \xi^2 + \eta^2)^{1/2} \mathcal{F}g.$ 
\end{lemma}

The previous result follows as a consequence of preliminary estimates involving only the leading part of the operator $P,$ namely $\partial_t + \partial_x^3 + \partial_x \partial_y^2.$ These are contained in the following lemma.
\begin{lemma}
\label{two_estimates}
	Let $w$ be as in Lemma~\ref{lemma:persistence} and let $D:=\R^2 \times [0,1].$ Then
	\begin{enumerate}
		\item For $\lambda>0$,
		\begin{multline}
		\label{eq:first}
			||e^{\lambda |x|}w||_{L^\infty_t L^2_{xy}(D)}
		\leq 
		||e^{\lambda |x|} w(0))||_{L^2(\R^2)} 
		+ ||e^{\lambda |x|} w(1)||_{L^2(\R^2)}\\
		+ ||e^{\lambda |x|} (\partial_t + \partial_x^3 + \partial_y\partial_x^2) w||_{L^1_t L^2_{xy}(D)}.
		\end{multline}
		\item There exists $c>0,$ independent of $K,$ such that for $\lambda\geq 1,$
		\begin{multline}
		\label{eq:second}
			||e^{\lambda |x|}L w||_{L^\infty_x L^2_{y t} (D)}
		\leq 
		c \lambda^2 \left(
		||J^3(e^{\lambda |x|} w(0)) ||_{L^2(\R^2)} 
		+ ||J^3(e^{\lambda |x|} w(1)) ||_{L^2(\R^2)}
		\right)\\
		+ c ||e^{\lambda |x|} (\partial_t + \partial_x^3 + \partial_y\partial_x^2) w ||_{L^1_x L^2_{ty}(D)},
		\end{multline}
		where $L$ denotes any operator in the set $\{\partial_x, \partial_y, \partial_x^2, \partial_{x y}, \partial_y^2\}.$
	\end{enumerate}
\end{lemma}  
Before sketching the proof of Lemma~\ref{two_estimates} we want to show how to get Lemma~\ref{lemma:persistence} once Lemma~\ref{two_estimates} is available.
\begin{proof}[of Lemma~\ref{lemma:persistence}]
	From estimate~\eqref{eq:first}, using the immediate estimates $||\cdot ||_{L^2(D)}\leq || \cdot ||_{L^\infty_t L^2_{x y}(D)}$ and $|| \cdot ||_{L^1_t L^2_{x y}(D)}\leq ||\cdot ||_{L^2(D)},$ adding and subtracting the lower order part of the operator $P$ and making also used of the H\"older inequality, one has
	\begin{multline}
		||e^{\lambda |x|}w||_{L^2(D)}
		\leq 
		||e^{\lambda |x|} w(0))||_{L^2(\R^2)} 
		+ ||e^{\lambda |x|} w(1)||_{L^2(\R^2)}\\
		+ ||e^{\lambda |x|} (\partial_t + \partial_x^3 + \partial_y\partial_x^2 + a_1\partial_x + a_0) w||_{L^1_t L^2_{xy}(D)}\\
		+ ||a_1 ||_{L^2_x L^\infty_{ty}(D)} || e^{\lambda |x|} \partial_x w||_{L^\infty_x L^2_{ty}(D)}\\
		+ ||a_0 ||_{L^\infty(D)} || e^{\lambda |x|} w||_{L^2(D)}.
	\end{multline}
	Similarly, from~\eqref{eq:second} one gets
	\begin{multline}
		||e^{\lambda |x|}L w||_{L^\infty_x L^2_{y t} (D)}
		\leq 
		c \lambda^2 \left(
		||J^3(e^{\lambda |x|} w(0)) ||_{L^2(\R^2)} 
		+ ||J^3(e^{\lambda |x|} w(1)) ||_{L^2(\R^2)}
		\right)\\
		+ c ||e^{\lambda |x|} (\partial_t + \partial_x^3 + \partial_y\partial_x^2 + a_1\partial_x + a_0) w ||_{L^1_x L^2_{ty}(D)}\\
		+ c||a_1 ||_{L^1_x L^\infty_{ty}(D)} || e^{\lambda |x|} \partial_x w||_{L^\infty_x L^2_{ty}(D)}\\
		+ c||a_0 ||_{L^2_x L^\infty_{ty}(D)} || e^{\lambda |x|} w||_{L^2(D)}.
	\end{multline}
	Summing the two previous estimates together we obtain
	\begin{multline}
		||e^{\lambda |x|}w||_{L^2(D)}
		+ \sum_{0<k+l\leq 2}||e^{\lambda |x|}\partial_x^k \partial_y^l w||_{L^\infty_x L^2_{y t} (D)}\\
		\leq c \lambda^2 \left(
		||J^3(e^{\lambda |x|} w(0)) ||_{L^2(\R^2)} 
		+ ||J^3(e^{\lambda |x|} w(1)) ||_{L^2(\R^2)}
		\right)\\
		+ c ||e^{\lambda |x|} (\partial_t + \partial_x^3 + \partial_y\partial_x^2 + a_1\partial_x + a_0) w ||_{L^1_t L^2_{xy}(D)\cap L^1_x L^2_{ty}(D)}\\
		+c ||a_1 ||_{L^2_x L^\infty_{ty}(D)\cap L^1_x L^\infty_{ty}(D)} || e^{\lambda |x|} \partial_x w||_{L^\infty_x L^2_{ty}(D)}\\
		+ c||a_0 ||_{L^\infty(D)\cap L^2_x L^\infty_{ty}(D)} || e^{\lambda |x|} w||_{L^2(D)}.
	\end{multline}
	Finally, estimate~\eqref{eq:persistence} follows absorbing the last two term in the right-hand side of the previous inequality into the corresponding terms in the left-hand side making use of the smallness of the norms of $a_0$ and $a_1$ in the corresponding spaces.
\end{proof}

Now we shall comment on the proof of Lemma~\ref{two_estimates}.
This result is related to the following lemma proved in~\cite{B_I_M} (Theorem 1.2 there).
\begin{lemma}
\label{lemma:persistenceBIM}
	Let $w$ be as in Lemma~\ref{lemma:persistence} and let $D:=\R^2 \times [0,1].$ Then
	\begin{enumerate}
		\item For $\lambda>0$ and $\beta>0$,
		\begin{multline}
		\label{eq:firstBIM}
			||e^{\lambda |x|} e^{\beta |y|} w||_{L^\infty_t L^2_{xy}(D)}
		\leq 
		||e^{\lambda |x|} e^{\beta |y|} w(0)||_{L^2(\R^2)} 
		+ ||e^{\lambda |x|} e^{\beta |y|} w(1)||_{L^2(\R^2)}\\
		+ ||e^{\lambda |x|} e^{\beta |y|} (\partial_t + \partial_x^3 + \partial_y\partial_x^2) w||_{L^1_t L^2_{xy}(D)}.
		\end{multline}
		\item There exists $c>0,$ independent of $K,$ such that for $\beta\geq 1$ and $\lambda\geq 7\beta$
		\begin{multline}
		\label{eq:secondBIM}
			||e^{\lambda |x|}e^{\beta |y|}L w||_{L^\infty_x L^2_{y t} (D)}\\
		\leq 
		c (\lambda^2 + \beta^2) \left(
		||J^3(e^{\lambda |x|} e^{\beta |y|} w(0)) ||_{L^2(\R^2)} 
		+ ||J^3(e^{\lambda |x|} e^{\beta |y|} w(1)) ||_{L^2(\R^2)}
		\right)\\
		+ c ||e^{\lambda |x|} e^{\beta |y|} (\partial_t + \partial_x^3 + \partial_y\partial_x^2) w ||_{L^1_x L^2_{ty}(D)},
		\end{multline}
		where $L$ denotes any operator in the set $\{\partial_x, \partial_y, \partial_x^2, \partial_{x y}, \partial_y^2\}.$
	\end{enumerate}
\end{lemma}

The conditions about $\beta$ in the previous result seem to exclude the case $\beta=0$ which corresponds to the estimates in Lemma~\ref{two_estimates} we aim at proving.  Nevertheless a closer look into the proof of Lemma~\ref{lemma:persistenceBIM} in~\cite{B_I_M} shows that the restrictions on $\beta$ are needed only in order to deal with  some technical difficulties which arise once one works with the \emph{full} operator $e^{\lambda |x|} e^{\beta |y|} (\partial_t + \partial_x^3 + \partial_y\partial_x^2).$  On the contrary, if one is interested just in the operator $e^{\lambda |x|}  (\partial_t + \partial_x^3 + \partial_y\partial_x^2),$ as we are, and mimics the proof of Lemma~\ref{lemma:persistenceBIM} in this situation, it turns out that the mentioned pathologies do not show up and  Lemma~\ref{two_estimates} follows in fact more easily than Lemma~\ref{lemma:persistenceBIM}. To see that we start with some preliminary observations. If we introduce the operator
\begin{equation*}
	\begin{split}
	H_{\lambda, \beta}:&=e^{\lambda x} e^{\beta y}(\partial_t + \partial_x^3 + \partial_x \partial_y^2) e^{-\lambda x} e^{-\beta y}\\
	&=\partial_t + (\partial_x -\lambda)^3 + (\partial_x -\lambda)(\partial_y -\beta)^2,
	\end{split}
\end{equation*}
and the quantity 
\begin{equation*}
	h=H_{\lambda, \beta} e^{\lambda x} e^{\beta y} w,
\end{equation*}
one realizes that the proof of the estimates in Lemma~\ref{lemma:persistenceBIM} is related to the two next results on the boundedness of the inverse of $H_{\lambda, \beta}.$ More precisely~\eqref{eq:firstBIM} follows from
\begin{equation}
\label{eq:boundedness}
||H_{\lambda, \beta}^{-1} h||_{L^\infty_t L^2_{xy}(\R^3)}
		\leq ||h||_{L^1_t L^2_{xy}(\R^3)}
\end{equation} 
and~\eqref{eq:secondBIM} (choosing, for instance, $L=\partial_x^2$ for notational simplicity) follows from 
\begin{equation}
\label{eq:boundednessL}
	||(\partial_x-\lambda)^2 H_{\lambda, \beta}^{-1} h||_{L^\infty_x L^2_{t y} (\R^3)}
		\leq 
	 c ||h||_{L^1_x L^2_{ty}(\R^3)},
\end{equation}
where the last implication is better understood observing that
\begin{equation*}
	\begin{split}
	e^{\lambda x} e^{\beta y} \partial_x^2 w
	&= e^{\lambda x} e^{\beta y} \partial_x^2 e^{-\lambda x} e^{-\beta y} e^{\lambda x} e^{\beta y} w
	=(\partial_x-\lambda)^2 e^{\lambda x} e^{\beta y} w\\
	&=(\partial_x-\lambda)^2  H_{\lambda, \beta}^{-1} h.
	\end{split}
\end{equation*}
For the proof of~\eqref{eq:boundedness} and~\eqref{eq:boundednessL} one studies the symbol of the operators $H_{\lambda, \beta}^{-1}$ and $(\partial_x - \lambda)^2 H_{\lambda, \beta}^{-1},$ namely   
\begin{equation*}
	m_0(\xi, \eta, \tau):=\frac{1}{i\tau + (i\xi-\lambda)^3 + (i \xi-\lambda)(i\eta -\beta)^2}
\end{equation*}
and 
\begin{equation*}
	m_2(\xi, \eta, \tau):=(i\xi-\lambda)^2 m_0(\xi, \eta,\tau),
\end{equation*}
respectively and the main ingredient in the proof is basically the explicit expression of the inverse Fourier transform of the following function
\begin{equation}
\label{gtau}
	g(\tau)= \frac{-i}{\tau + a + i b}, \quad b\neq 0,
\end{equation}
which is
\begin{equation}
\label{explicitF}
	\mathcal{F}^{-1}g(t)=
	\begin{system}
		&\sqrt{2\pi}\chi_{(0,+\infty)}(t) e^{-ita} e^{tb} \quad &\text{if}\, b<0,\\
		&-\sqrt{2\pi}\chi_{(-\infty,0)}(t) e^{-ita} e^{tb} &\text{if}\, b>0.
	\end{system}
\end{equation}
Notice that the size of the expression above is bounded by $\sqrt{2\pi}.$
For the proof of~\eqref{eq:boundedness} one starts observing that the symbol $m_0(\xi, \eta, \tau)$ has the following favorable form: 
\begin{equation}
\label{m0_dec}
m_0(\xi, \eta, \tau)=\frac{-i}{\tau + a(\xi, \eta) + i b(\xi, \eta)}, 
\end{equation}
which is the same as $g$ in~\eqref{gtau} and where
\begin{equation*}
\begin{split}
&a(\xi, \eta)=-\xi^3 + (3\lambda^2 + \beta^2)\xi -\xi\eta^2 + 2\lambda \beta \eta,
\\
&b(\xi, \eta)=-3\lambda \xi^2 -2\beta \xi \eta -\lambda \eta^2 + \lambda^3 + \lambda \beta^2.
\end{split}
\end{equation*}

From this fact, applying Plancherel's formula, using that 
\begin{equation*}
	\mathcal{F}^{-1}[f g](t)= \frac{\mathcal{F}^{-1}[f](t) \ast \mathcal{F}^{-1}[g](t)}{\sqrt{2\pi}},
\end{equation*}
 using~\eqref{explicitF} and, finally, the Minkowski integral inequality, one gets
\begin{equation*}
\begin{split}
||H_{\lambda, \beta}^{-1}h(t)||_{L^2_{x,y}(\R^2)}
&=||\mathcal{F}^{-1}_\tau [m_0(\xi,\eta, \cdot_\tau) \mathcal{F}h(\xi,\eta,\cdot_\tau)](t)||_{L^2_{\xi,\eta}(\R^2)}\\
&=\frac{1}{\sqrt{2\pi}}||\mathcal{F}^{-1}_\tau [m_0(\xi, \eta, \cdot_\tau)](t)\ast_t \mathcal{F}_{x,y}[h(\cdot_x, \cdot_y, t)](\xi,\eta)||_{L^2_{\xi,\eta}(\R^2)}\\
&\leq \int_{\R} || \mathcal{F}_{xy}h(\cdot_\xi, \cdot_\eta, s) ||_{L^2_{\xi,\eta}(\R^2)}\, ds\\
&= ||h||_{L^1_t L^2_{x, y}(\R^3)},
\end{split}
\end{equation*}
which gives~\eqref{eq:boundedness}.

We stress that ensuring that the measure of the set $\{(\xi, \eta)\colon b(\xi, \eta)=0\}$ is zero in $\R^2,$ which is needed in order to use meaningfully~\eqref{explicitF} in the previous estimate, requires that the parameters $\lambda$ and $\beta$ are not \emph{both} equal to zero, then if one fixes $\beta=0$ from the beginning, it is enough asking $\lambda\neq 0$ to make the argument work.  This shows that the proof of estimate~\eqref{eq:boundedness} (and so of~\eqref{eq:firstBIM} since it is its consequence) can be extended to the case $\beta=0$ if one assumes, for instance, $\lambda>0.$ This gives our desired estimate~\eqref{eq:first} in Lemma~\ref{two_estimates}.   

The proof of~\eqref{eq:secondBIM} is more delicate and needs to perform a preliminary partial fraction decomposition in order to recognize that $m_2$ can be written again as $g$ in~\eqref{gtau} and then to be able to use~\eqref{explicitF}. First, for our convenience, we write
\begin{equation}
\label{fraction}
	m_2(\xi,\eta,\tau)= \frac{P(v)}{Q(v)},
\end{equation}
where $P(v)=-iv^2$ and $Q(v)=v^3 + w^2 v-\tau,$ with 
\begin{equation*}
	v:= \xi +i \lambda
	\quad \text{and} \quad
	w:=\eta + i \beta.
\end{equation*}
\emph{If} $w$ and $\tau$ are such that the polynomial \emph{$Q$ has no multiple roots}, one can obtain the estimate~\eqref{eq:secondBIM} basically in the same way as in the previous case. Indeed in this case one can apply a partial fraction decomposition and gets
\begin{equation}
\label{fraction_dec}
	\frac{P(v)}{Q(v)}=\sum_{j=1}^3 \frac{P(v_j)}{Q'(v_j)(v-v_j)}
	=\sum_{j=1}^3 A_j \frac{-i}{\xi + a_j(\eta, \tau) + i b_j(\eta, \tau)},
\end{equation}
where $v_j,$ $j=1,2,3,$ are the different roots of 
$Q,$
\begin{equation*} 
	A_j:=\frac{-i v_j^2}{Q'(v_j)},
\end{equation*}
moreover $a_j(\eta, \tau)=-\Re(v_j)$ and $b_j(\eta,\tau)=\lambda - \Im(v_j).$
Combining~\eqref{fraction} together with~\eqref{fraction_dec} one obtains
\begin{equation*}
	m_2(\xi,\eta,\tau)=\sum_{j=1}^3 A_j \frac{-i}{\xi + a_j(\eta, \tau) + i b_j(\eta, \tau)}.
\end{equation*}
Here, we clearly recognize the same structure as the one of $m_0$ in the previous case, namely~\eqref{m0_dec}, and then, as expected, the proof of~\eqref{eq:boundednessL} follows analogously. For this reason we do not comment further on this. 

\emph{Nevertheless} one should notice that \emph{$Q(v)$ does have multiple roots}. More precisely, $Q(v)$ has multiple roots if and only if 
\begin{equation}
\label{cond_mr}
	\frac{\tau^2}{w^6}=-\frac{4}{27}.
\end{equation}
If this happens, a more careful analysis is needed in order to get~\eqref{eq:secondBIM} and it is here that the restrictions $\lambda\geq 1$ and $\beta\geq 7\lambda$ that appear in the statement of~\eqref{eq:secondBIM} in Lemma~\ref{lemma:persistenceBIM} are required. 

On the other hand, we stress that for our purpose of proving estimate~\eqref{eq:second}, we need to impose $\beta=0.$ Notice that setting $\beta=0$   
rules out the possibility for $Q$ to have multiple roots, indeed if $\beta=0,$ then $w$ is a purely real number and condition~\eqref{cond_mr} is never satisfied. This means that in this situation the proof follows without any other parameter restrictions. This concludes the sketch of the proof of Lemma~\ref{two_estimates}.

\begin{proof}(of Theorem~\ref{upper_bound})
	We consider a $C^\infty(\R)$ truncation function $\mu_R$ such that $\mu_R(x)=0$ if $|x|< R$ and $\mu_R(x)=1$ if $|x|> (18R-1)/4.$
	
	We define 
	\begin{equation*}
		w(x,y,t):=\mu_R(x) v(x,y,t), \qquad (x,y)\in \R^2, t\in[0,1],
	\end{equation*}
	where $v$ satisfies~\eqref{ZKv_general}.
	
	It is easy to see that $w$ satisfies
	\begin{equation*}
		(\partial_t + \partial_x^3 + \partial_x\partial_y^2 + a_1 \partial_x + a_0)w=e_R,
	\end{equation*}
	where
	\begin{equation*}
		e_R=\partial_x^3 \mu_R v + 3\partial_x^2 \mu_R \partial_x v+3 \partial_x \mu_R \partial_x^2 v + \partial_x \mu_R \partial_y^2 v + a_1 \partial_x \mu_R v.
	\end{equation*}
	
	The next step would be to apply estimate~\eqref{eq:persistence} of Lemma~\ref{lemma:persistence} to $w.$ Observe that in order to do that $a_0$ and $a_1$ are required to have small norms in the corresponding spaces, but, as a matter of fact, this is not necessarily true under the generous hypotheses of Theorem~\ref{upper_bound}. For this reason we introduce an auxiliary function  $\widetilde{\mu}_R$ such that $\widetilde{\mu}_R \mu_R=\mu_R$ and $\widetilde{a}_j:=a_j \widetilde{\mu}_R,$ for $j=0,1,$ have small norms for $R$ sufficiently large.
	
	It is easy to show that
	\begin{equation*}
		(\partial_t + \partial_x^3 + \partial_x\partial_y^2 + \widetilde{a}_1 \partial_x + \widetilde{a}_0)w=\widetilde{e}_R,
	\end{equation*}
	where $\widetilde{e}_R$ has the same form as $e_R$ but with $a_j$ replaced by $\widetilde{a}_j,$ for $j=0,1.$
	
	Now we can apply estimate~\eqref{eq:persistence} of Lemma~\ref{lemma:persistence} to $w.$ This gives
	\begin{multline}
	\label{preliminary}
		||e^{\lambda |x|}w||_{L^2(D)}
		+\sum_{0<k+l\leq 2}||e^{\lambda |x|}\partial_x^k \partial_y^l w||_{L^\infty_x L^2_{y t} (D)}\\
		\leq 
		c \lambda^2 \left(
		||J^3(e^{\lambda |x|} w(0)) ||_{L^2(\R^2)} 
		+ ||J^3(e^{\lambda |x|} w(1)) ||_{L^2(\R^2)}
		\right)\\
		+ c ||e^{\lambda |x|} \widetilde{e}_R ||_{L^1_t L^2_{xy}(D)\cap L^1_x L^2_{ty}(D)}.
	\end{multline}
	
	We want to estimate the right-hand side of the previous inequality. Let us consider the term $||J^3(e^{\lambda |x|} w(0)) ||_{L^2(\R^2)}$ first. 
	The following interpolation result proved in~\cite{B_I_M0} (Lemma 1 there) will be useful for our purpose.
	\begin{lemma}
	\label{interpolation}
		For $s>0$ and $\beta>0$ let $f\in H^s(\R^2)\cap L^2(e^{2\beta |x|}dx dy).$ Then for $\theta\in [0,1]$ one has
		\begin{equation}
		\label{eq:interpolation}
			|| J^{\theta s}(e^{(1-\theta)\beta |x|} f) ||_{L^2(\R^2)}
			\leq c ||J^s f||_{L^2(\R^2)}^\theta ||e^{\beta |x|} f||_{L^2(\R^2)}^{1-\theta},
		\end{equation}
		with a constant $c>0$ depending on $s$ and $\beta.$
	\end{lemma}
	We want to apply the previous estimate to $f=w(0)$ in the case $s=4$ and setting $\beta=4\lambda.$ In order to do that we shall show that $w(0)$ satisfies the regularity hypotheses of Lemma~\ref{interpolation}. 
	
	First of all we show that $w(0)\in H^4(\R^2).$ Notice that since $\mu_R$ and its derivatives are bounded by a constant independent of $R,$ we easily have
	\begin{equation*}
		||J^4 w(0)||_{L^2(\R^2)}\leq c ||J^4 v(0)||_{L^2(\R^2)}. 
	\end{equation*}
	Moreover, by hypothesis $v\in C([0,1]; H^4(\R^2)),$ which ensures the finiteness of the last norm. Thus $w(0)\in H^4(\R^2).$ 
	
	Now it is left to show that $||e^{4\lambda |x|} w(0) ||_{L^2(\R^2)}$ is finite. Using that $w$ is supported in the set $\{(x,y,t)\colon |x|\geq R, t\in [0,1]\}$ and again that $\mu_R$ is bounded by a constant independent of $R$ one has
	\begin{equation}
	\label{est_w(0)}
		||e^{4\lambda |x|} w(0) ||_{L^2(\R^2)}=
		||e^{4\lambda |x|} w(0) ||_{L^2(|x|\geq R)}
		\leq c ||e^{4\lambda |x|} v(0) ||_{L^2(|x|\geq R)}.
	\end{equation}
	Choosing 
	\begin{equation}
	\label{lambda}
		\lambda= \frac{2aR^{3/2}}{18R-1},
	\end{equation}
	for $R$ sufficiently large, one has 
	\begin{equation*}
		e^{4\lambda |x|}= e^{4 \frac{2aR^{3/2}}{18R-1} |x|}\leq e^{\frac{a}{2}R^{1/2} |x|}\leq e^{\frac{a}{2}|x|^{3/2}} \qquad \text{for}\, |x|\geq R.
	\end{equation*}
	This fact implies the bound
	\begin{equation*}
	 ||e^{4\lambda |x|} v(0) ||_{L^2(|x|\geq R)}\leq ||e^{\frac{a}{2} |x|^{3/2}} v(0)||_{L^2(|x|\geq R)}.
	\end{equation*}
	Since, by assumption, $v(0)\in L^2(e^{a |x|^{3/2}} dx dy),$ then $||e^{\frac{a}{2} |x|^{3/2}} v(0)||_{L^2(|x|\geq R)}$ is finite and so is $||e^{4\lambda |x|} w(0) ||_{L^2(\R^2)}$ due to~\eqref{est_w(0)}.
	
	Now we are in position to apply Lemma~\ref{interpolation} to estimate $||J^3(e^{\lambda |x|} w(0)) ||_{L^2(\R^2)}.$ More precisely, choosing $f=w(0),$ $s=4,$ $\beta=4\lambda$ and $\theta =3/4$ in~\eqref{eq:interpolation}, we have
	\begin{equation}
	\label{bound_w0}
		\begin{split}
		||J^3(e^{\lambda |x|} w(0)) ||_{L^2(\R^2)}&\leq 
		c ||J^4 w(0)||_{L^2(\R^2)}^{3/4} ||e^{4\lambda |x|} w(0)||_{L^2(\R^2)}^{1/4}\\
		&\leq c ||J^4 v(0)||_{L^2(\R^2)}^{3/4} ||e^{\frac{a}{2} |x|^{3/2}} v(0)||_{L^2(|x|\geq R)}\\
		&\leq C.
		\end{split}
	\end{equation}
	A similar argument shows that 
	\begin{equation}
	\label{bound_w1}
		||J^3(e^{\lambda |x|} w(1)) ||_{L^2(\R^2)}\leq C.
	\end{equation}
	
	It remains to bound the term $||e^{\lambda |x|} \widetilde{e}_R||_{L^1_t L^2_{xy}(D)\cap L^1_x L^2_{ty}(D)}$ in~\eqref{preliminary}.
	
	Since $\widetilde{e}_R$ is supported in $\Omega_R:=\{(x,y,t) \colon R\leq |x|\leq (18R-1)/4, y\in \R, t\in [0,1]\},$ we find that
	\begin{equation}
	\label{bound_eR}
		\begin{split}
		||e^{\lambda |x|} \widetilde{e}_R||&_{L^1_t L^2_{xy}(D)\cap L^1_x L^2_{ty}(D)}\\
		&\leq e^{\lambda \frac{18R-1}{4}} || \widetilde{e}_R \chi_{\Omega_R}||_{L^1_t L^2_{xy}(D)\cap L^1_x L^2_{ty}(D)}\\
		&\leq c e^{\lambda \frac{18R-1}{4}} 
		|| (|v| + |\partial_x v| + |\partial_y v| + |\Delta v|) \chi_{\Omega_R}||_{L^1_t L^2_{xy}(D)\cap L^1_x L^2_{ty}(D)}\\
		&\leq c R^{1/2} e^{\lambda \frac{18R-1}{4}},
		\end{split}
	\end{equation}
	where in the last inequality we have applied the H\"older inequality. Here and in the following we use the same letters, $c$ or $C$, to denote diverse constants.
	
	Using~\eqref{bound_w0},~\eqref{bound_w1} and~\eqref{bound_eR} in~\eqref{preliminary} we have
	\begin{equation}
	\label{last}
		\begin{split}
		||e^{\lambda |x|}w||_{L^2(D)}
		+\sum_{0<k+l\leq 2}||e^{\lambda |x|}\partial_x^k \partial_y^l w||_{L^\infty_x L^2_{y t} (D)}
		&\leq 
		C\lambda^2 + C R^{1/2} e^{\lambda \frac{18 R-1}{4}}\\ 
		&\leq
		C\lambda^2 R^{1/2} e^{\lambda \frac{18 R-1}{4}}.
		\end{split}
	\end{equation}
	Let $Q_{18R}:=\{(x,y)\colon 18R - 1\leq |x|\leq 18 R \wedge |y|\leq 2\},$ for $R$ sufficiently large, one has $Q_{18R} \subseteq \{|x|\geq R, y\in \R\},$ where the last one is the set in which $w(t)$ is supported, for any $t\in [0,1].$ Moreover in $Q_{18R}\times [0,1]$ we have $w=v.$ From these facts, together with~\eqref{last} and the H\"older inequality, we get
	\begin{equation}
	\label{preliminary2}
		\begin{split}
			||e^{\lambda |x|}v||&_{L^2(Q_{18R}\times [0,1])}
		+\sum_{0<k+l\leq 2}||e^{\lambda |x|}\partial_x^k \partial_y^l v||_{L^2(Q_{18R}\times [0,1])}\\
		&\leq R^{1/2}\Big(||e^{\lambda |x|}w||_{L^2(D)}
		+\sum_{0<k+l\leq 2}||e^{\lambda |x|}\partial_x^k \partial_y^l w||_{L^\infty_x L^2_{y t} (D)}\Big)
		\\ 
		&\leq C \lambda^2 R e^{\lambda \frac{18 R-1}{4}}.
		\end{split}
	\end{equation}
	Notice that, if $|x|\geq 18R-1,$ one has
	\begin{equation}
	\label{left}
		\lambda |x|\geq \lambda(18R-1)=2a R^{3/2},
	\end{equation}
	where in the last identity we have used our explicit choice of $\lambda,$ namely~\eqref{lambda}.
	
	Moreover, again using~\eqref{lambda} and that, for sufficiently large $R,$ $\lambda\geq 1,$ we get
	\begin{equation}
	\label{right}
		\lambda^2 R e^{\lambda \frac{18R-1}{4}}
		\leq C a^2 R^2 e^{\lambda \frac{18R-1}{4}}
		\leq C_a e^{(\lambda + 1) \frac{18R-1}{4}}
		\leq C_a e^{2\lambda \frac{18R-1}{4}}
		= C_a e^{a R^{3/2}}
	\end{equation}
	Plugging~\eqref{left} and~\eqref{right} in~\eqref{preliminary2} implies
	\begin{equation*}
		e^{2aR^{3/2}}
		\Big(||v||_{L^2(Q_{18R}\times [0,1])}
		+\sum_{0<k+l\leq 2}||\partial_x^k \partial_y^l v||_{L^2(Q_{18R}\times [0,1])}\Big)
		\leq C_a e^{a R^{3/2}},
	\end{equation*}
	or, equivalently,
	\begin{equation*}
		||v||_{L^2(Q_{18R}\times [0,1])}
		+\sum_{0<k+l\leq 2}||\partial_x^k \partial_y^l v||_{L^2(Q_{18R}\times [0,1])} \leq C_a e^{-a R^{3/2}}.
	\end{equation*}
	A rescale in $R$ gives
	\begin{equation*}
		||v||_{L^2(Q_{R}\times [0,1])}
		+\sum_{0<k+l\leq 2}||\partial_x^k \partial_y^l v||_{L^2(Q_{R}\times [0,1])} \leq C_a e^{-a \big(\frac{R}{18}\big)^{3/2}}.
	\end{equation*}
	The desired upper bound~\eqref{eq:upper_bound} follows once one notices that the left-hand side of the previous expression can be bounded from below by $A_R(v).$ This completes the proof. 
	\end{proof}

\begin{acknowledgement}
The author gratefully acknowledges financial support by the Deutsche Forschungsgemeinschaft (DFG) through CRC 1173. The author would like to thank also L. Fanelli and F. Linares for valuable comments and suggestions on the draft of the manuscript. 
\end{acknowledgement}
\section*{Appendix}
\addcontentsline{toc}{section}{Appendix}

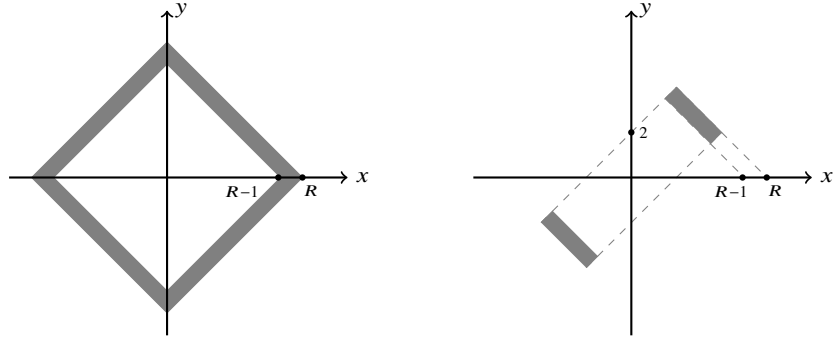
\begin{figure}[ht]
\centering
\begin{tikzpicture}[scale=.3]
\filldraw[gray] (0,5)-- (0,6)--(6,0)--(5,0);
\filldraw[gray] (0,-5)-- (0,-6)--(6,0)--(5,0);
\filldraw[gray] (0,-5)-- (0,-6)--(-6,0)--(-5,0);
\filldraw[gray] (0,5)-- (0,6)--(-6,0)--(-5,0);
\coordinate [label=below left: \textcolor{black}{$\scriptstyle R-1$}] (a) at (4.4,0); 
\fill[black] (4.93,0) circle (4pt);
\coordinate [label=below right: \textcolor{black}{$\scriptstyle R$}] (b) at (5.7,0); 
\fill[black] (6,0) circle (4pt);
\draw[thick,->](-7,0)--(8,0) node[right]{$x$};
\draw[thick,->](0,-7)--(0,7.4) node[right]{$y$};
\end{tikzpicture}
\qquad \qquad
\begin{tikzpicture}[scale=.3]
\draw[domain=-4:2,color=gray,dashed,samples=2]
plot(\x, {\x +2});
\draw[domain=-2:4,color=gray,dashed,samples=2]
plot(\x, {\x -2});
\draw[domain=1.5:5,color=gray,dashed,samples=2]
plot(\x, {-\x +5});
\draw[domain=2:6,color=gray,dashed,samples=2]
plot(\x, {-\x +6});

\filldraw[gray] (1.5,3.5)-- (2,4)--(4,2)--(3.5,1.5);
\filldraw[gray] (-4,-2)-- (-3.5,-1.5)--(-1.5,-3.5)--(-2,-4);
\coordinate [label=below: \textcolor{black}{$\scriptstyle R-1$}] (a) at (4.4,0); 
\fill[black] (4.93,0) circle (4pt);
\coordinate [label=below right: \textcolor{black}{$\scriptstyle R$}] (b) at (5.7,0); 
\fill[black] (6,0) circle (4pt);
\coordinate [label=right: \textcolor{black}{$\scriptstyle 2$}] (b) at (0,2); 
\fill[black] (0,2) circle (4pt);
\draw[thick,->](-7,0)--(8,0) node[right]{$x$};
\draw[thick,->](0,-7)--(0,7.4) node[right]{$y$};
\end{tikzpicture}
\caption{The regions $Q_R$ (left) and $\widetilde{Q}_R$ (right)}
\label{QRold} 
\end{figure}

This appendix is concerned with fixing an inaccuracy we spotted in the proof of Theorem 1.4 in~\cite{C_F_L} (Theorem~\ref{thm:intermediate} here). 
As already mentioned, this result follows from the comparison between a lower bound and an upper bound for the following quantity
\begin{equation}\label{eq:quantity}
	A_R(v):=\Big( \int_0^1 \int_{Q_R} (|v|^2 + |\nabla v|^2 + |\Delta v|^2)\, dx dy dt \Big)^{1/2}, 
\end{equation}
where $v:=u_1-u_2$ is the difference of two solutions $u_1$ and $u_2$ of~\eqref{ZKsymmetric} and where $Q_R:=\{(x,y)\colon |x+y|\leq R \wedge |x-y|\leq R\}\setminus \{(x,y)\colon |x+y|<R-1 \wedge |x-y|<R-1\}$ (see Figure~\ref{QRold} left).  It turned out that this definition of $Q_R$ was not totally correct. Nevertheless the proof in~\cite{C_F_L} remains valid once one implements the following few modifications:
\begin{itemize}
	\item We replace the domain of integration $Q_R$ in~\eqref{eq:quantity} by the alternative domain $\widetilde{Q}_R:=\{(x,y)\colon R-1\leq |x+y|\leq R \wedge |x-y|\leq 2\}$  (see Figure~\ref{QRold} right). 
	In passing, observe that a suitable rotation of the Cartesian axes brings $\widetilde{Q}_R$ into the corresponding domain, sketched in Figure~\ref{fig:an_dom}, that we introduced in the previous sections for proving Theorem~\ref{main} without passing through the symmetrization procedure. This fact, according to the change of variables~\eqref{change_variable}, makes the new domain choice $\widetilde{Q}_R$ very natural.
	\item In the proof of the lower bound (Theorem 2.1 in~\cite{C_F_L}) the authors introduced the following auxiliary modified-$v$ function
\begin{equation*}
	g(x,y,t)=\theta_R(x,y) \mu\left(\frac{x}{R} + \varphi(t), \frac{y}{R} + \varphi(t) \right) v(x,y,t),
\end{equation*}
where $\theta_R\in C^\infty_0(\R^2)$ is chosen such that $\theta_R(x,y)=1$ on $\{(x,y)\colon |x+y|< R-1 \wedge |x-y|<R-1\}$ and $\theta_R(x,y)=0$ on $\R^2\setminus \{(x,y)\colon |x+y|\leq R \wedge |x-y|\leq R\}$ (this definitions corrects also a typo in~\cite{C_F_L}), $\mu\in C^\infty(\R^2)$ is such that $\mu(x,y)=0$ if $\sqrt{x^2+y^2}<1$ and $\mu(x,y)=1$ if $\sqrt{x^2+y^2}>2$ and $v:=u_1-u_2$ is again the difference of the two solutions $u_1$ and  $u_2$ of~\eqref{ZKsymmetric}.  In order to be consistent with the new choice $\widetilde{Q}_R$ of the domain and with the decay assumption~\eqref{decay_symmetric}, in analogy with the proof presented in the previous sections (see definition~\eqref{g}), we replace $\theta_R$ in the definition of $g$ with the function $\theta_{R,\varepsilon}\in C^\infty_0(\R^2)$ defined such that $\theta_{R, \varepsilon}(x,y)=1$ on $\{(x,y)\colon |x+y|< R-1 \wedge |x-y|<2\}$ and $\theta_{R,\varepsilon}(x,y)=0$ on $\R^2\setminus \{(x,y)\colon |x+y|\leq R \wedge |x-y|\leq 2+\varepsilon\}.$    
\end{itemize}
The proof of Theorem 1.4 in~\cite{C_F_L} with these modifications follows analogously to the proof of Theorem~\ref{main} here, therefore we will omit the details.

\begin{thebibliography}{99.}%
%
%
%
%
%


\bibitem{B_I_M0} E. Bustamante, P. Isaza, J. Mej\'ia, 
\emph{On the support of solutions to the Zakharov-Kuznetsov equation},
J. Differential Equations \textbf{251} (2011), no. 10,  pp. 2728-2736

%
\bibitem{B_I_M} E. Bustamante, P. Isaza, J. Mej\'ia, 
\emph{On uniqueness properties of solutions of the Zakharov-Kuznetsov equation},
J. Funct. Anal. \textbf{264} (2013), no. 11, pp. 2529-2549
%
\bibitem{B_J_M} E. Bustamante, J. Jim\'enez Urrea, J. Mej\'ia, 
\emph{On the unique continuation property of solutions of the three-dimensional Zakharov-Kuznetsov equation}
Nonlinear Anal. Real World Appl. \textbf{39} (2018), pp. 537-553
%
\bibitem{C_F_L} L. Cossetti, L. Fanelli, F. Linares, 
\emph{Uniqueness results for Zakharov-Kuznetsov equation}
Comm. Partial Differential Equations \textbf{44} (2019), no. 6, pp. 504-544
%
\bibitem{E_K_P_V} L. Escauriaza, C. E. Kenig, G. Ponce, L. Vega,
\emph{On uniqueness properties of solutions of the k-generalized KdV equations}, 
J. Funct. Anal \textbf{244} (2007), no. 2, pp. 504-535
%
\bibitem{F} A. V. Faminskii,
\emph{The Cauchy problem for the Zakharov-Kuznetsov equation},
Differential Equations \textbf{31} (1995), no. 6,  pp. 1002-1012
%
\bibitem{F_A} A. V. Faminskii, A. P. Antonova, 
\emph{On internal regularity of solutions to the initial value problem for the Zakharov-Kuznetsov equation},
in \textit{Progress in Partial Differential Equations} (Springer, Heidelberg, 2013), pp. 53-74
%
\bibitem{G_H} A. Gr\"unrock, S. Herr, 
\emph{The Fourier restriction norm method for the Zakharov-Kuznetsov equation},
Discrete Contin. Dyn. Syst. \textbf{34} (2014), no. 5, pp. 2061-2068
%
\bibitem{L_L_S}
	D. Lannes, F. Linares, J. C. Saut,
	\emph{The Cauchy problem for the Euler-Poisson system and derivation of the Zakharov-Kuznetsov equation},
	Progr. Nonlinear Differential Equations Appl. \textbf{84} (2013), pp 181-213
%
\bibitem{L_P}
	F. Linares,  A. Pastor,
	\emph{Well-posedness for the two-dimensional modified Zakharov-Kuznetsov equation},
	SIAM J. Math. Anal. \textbf{41} (2009), no. 4, pp. 1323-1339
%
\bibitem{L_PII}
	F. Linares,  A. Pastor,
	\emph{Local and global well-posedness for the 2D generalized Zakharov-Kuznetsov equation},
	J. Funct. Anal. \textbf{260} (2011), no. 4, pp. 1060-1085
%
\bibitem{L_P_S}
	F. Linares, A. Pastor, J. C. Saut,
	\emph{Well-posedness for the ZK equation in a cylinder and on the background of a KdV soliton},
	Comm. Partial Differential Equations \textbf{35} (2010), no. 9, pp. 1674-1689
%
\bibitem{M_P} L. Molinet, D. Pilod, 
\emph{Bilinear Strichartz estimates for the Zakharov-Kuznetsov equation and applications},
Ann. Inst. H. Poincaré Anal. Non Linéaire \textbf{32} (2010), no. 2, pp. 347-371
%
\bibitem{Z_K}
	V.E. Zakharov and E.A. Kuznetsov,
	\emph{On three-dimensional solitons},
	Sov. Phys. JETP 39 (1974), pp. 285-286.

%
\end{thebibliography}
\end{document}